%% file: main.tex
\documentclass[a4paper,twoside]{article}
\newcounter{the_style}
\input{main.macros}
\title{A Lower Bound on Arbitrary $f$--Divergences in Terms of the Total Variation}
\author{Jochen Br\"{o}cker\thanks{email: \texttt{broecker@pks.mpg.de}}\\
Max--Planck--Institut f\"{u}r Physik komplexer Systeme\\
N\"{o}thnitzer Strasse~34\\
01187 Dresden \\
Germany}
\begin{document}
\maketitle
\begin{abstract}
An important tool to quantify the likeness of two probability measures are $f$--divergences, which have seen widespread application in statistics and information theory.
An example is the total variation, which plays an exceptional role among the $f$--divergences.
It is shown that every $f$--divergence is bounded from below by a monotonous function of the total variation.
Under appropriate regularity conditions, this function is shown to be monotonous.
Remark: The proof of the main proposition is relatively easy, whence it is highly likely that the result is known. 
The author would be very grateful for any information regarding references or related work. 
\end{abstract}
\section{The total variation}
Let $(\Omega, \sigma)$ be a probability space.
A {\em signed measure } $\nu$ is a $\sigma$--additive set function with values in $\R \cup \{ -\infty, \infty\}$, and so that either $\nu > -\infty$ or $\nu < \infty$.
I will use the standard term {\em measure} if $\nu$ is nonnegative.
To any signed measure $\nu$, there corresponds a {\em Hahn--Jordan decomposition} of $\Omega$ into two measurable sets $P, N$ so that $P \cup N = \Omega$, $P \cap N = \emptyset$ and 
\beq{equ:10.10}
\nu^+(.) = \nu(. \cap P) , \qquad
\nu^-(.) = -\nu(. \cap N)
\eeq
are both (nonnegative) measures.
Obviously, $ \nu = \nu^+ - \nu^-$.
Furthermore, the representation
\beq{equ:10.20}
\nu^+(A) = \sup_{B \subset A} \nu(B) , \qquad
\nu^-(A) = -\inf_{B \subset A} \nu(B)
\eeq
holds for every measurable set $A$.
For a proof of these facts see~\cite{doob94}.
The measure  $ \vr{\nu} = \nu^+ + \nu^-$ is called the {\em variation measure} of $\nu$, which in turn defines the {\em total variation} $\| \nu \| = \vr{\nu}(\Omega)$.
If $ \nu(\Omega) = 0$, it follows easily from the previous statements that 
\beq{equ:10.30}
\vr{\nu}(\Omega) = 2 \sup_{B \in \sigma } |\nu(B)|.
\eeq
A {\em probability measure} is a measure $\mu$ so that $\mu(\Omega) = 1$.
For any two probability measures, $\mu$, $\nu$, the difference $\mu - \nu$ is a signed measure, and Equation~\eqref{equ:10.30} applies.
Hence,
\beq{equ:10.40}
\| \mu - \nu\| = \vr{\mu - \nu}(\Omega) = 2 \sup_{B \in \sigma } |\mu(B) - \nu(B)|.
\eeq
Obviously, $\| \mu - \nu\|$ is a metric for probability measures, namely the {\em total variation metric}, with Equation~\eqref{equ:10.40} providing two possible representations. 
If $\mu$ is absolutely continuous with respect to $\mu$, then there is a third representation, namely
\beq{equ:10.50}
\| \mu - \nu\| = \int |\odq{\mu}{\nu} - 1 | \dd \nu .
\eeq
\texttt{Proof of this fact}
\section{The $f$-divergences}
Equation~\eqref{equ:10.50} can be read as follows:
\beq{equ:20.10}
\| \mu - \nu \| = \int f(\odq{\mu}{\nu}) \dd \nu ,
\eeq
with $f(x) = | x - 1|$.
There is a way to generalise this approach by using other forms of $f$.
Let $f$ be a convex function on $\R_{\geq 0}$ that vanishes
at $x = 1$.
Let $\mu, \nu$ two probability measures with $\mu$ being absolutely continuous with respect to $\nu$ (which will be written as $\mu \ll \nu$).
The $f${\em{}--divergence} between $\mu$ and $\nu$ is given by 
\beq{equ:20.20}
D_f(\mu, \nu) = \int f(\odq{\mu}{\nu}) \dd \nu.
\eeq

For, if $\mu = \nu$ we have $\odq{\mu}{\nu} = 1$, we see that $f(\mu ,
\nu)$ vanishes in this case.
Furthermore, $D_f(\mu, \nu)$ is non-negative.
Indeed, by Jensen's inequality, 
\[
0 = f(1) 
=  f(\int \odq{\mu}{\nu} \dd \nu) 
\leq \int f(\odq{\mu}{\nu}) \dd \nu
= f(\mu , \nu).
\]
Note though that $f(\mu , \nu)$ may be infinite.
Furthermore $f(\mu , \nu)$ may vanish even if $\mu \neq \nu$.
To exclude this, further conditions on $f$ have to be imposed, for example as in the following
\begin{rmlist}[Lemma]
\label{rml:20.10}
Suppose there is an $a \in \R$ so that the function
\[
g(x) := f(x) - a(x - 1)
\]
is non-negative and vanishes only if $x = 1$, then $f(\mu , \nu)$ vanishes only if $\mu = \nu$.
\end{rmlist}
\begin{proof}
The function $g(x)$ is convex as well.
Furthermore $D_f(\mu, \nu) = D_g(\mu, \nu)$.
But since $g$ is non-negative,
\[
D_g(\mu , \nu) = \int g( \odq{\mu}{\nu}) \dd \nu 
\]
can only vanish if $g(\odq{\mu}{\nu})$ is identical to zero, which implies that 
$\odq{\mu}{\nu} = 1$ $\nu$-a.s.
But this means $\mu = \nu$.
\end{proof}
The concept of $f$-divergences was introduced by Csisz\'{a}r~\cite{csiszar67}, who also noted the result in Lemma~\ref{rml:20.10}. 
Common choices for $f$ are
\begin{align*}
(\sqrt{x} - 1)^2 & \qquad \mbox{Hellinger divergence~} \HE\\
|x - 1| & \qquad \mbox{total--variation divergence~} \TV \\
x \* \log(x) & \qquad \mbox{Kullback--Leibler divergence~} \KL \\
(x - 1)^2 & \qquad \mbox{Pearson divergence~} \PE
\end{align*}
The transformation $f^*(x) = x f(1 / x)$ yields a divergence $D_{f^*}$ which is equal to $D_{f}$ but with interchanged arguments.
Applying this transformation to the Kullback--Leibler divergence for example, we get a divergence which is also sometimes referred to as the Kullback--Leibler divergence, or alternatively as the Shannon divergence $\SH$.
The total variation divergence plays a central role, since all $f$--divergences allow for an estimate against $\TV$, as will be shown in the following proposition, which forms the main result of this short note.
\begin{rmlist}[Proposition]
\label{ftheorem}
For two probability measures $\mu, \nu$, it holds in general that
\[
f(1 + \frac{1}{2} \TV(\mu, \nu)) + f(1 - \frac{1}{2} \TV(\mu, \nu)) 
\leq D_f(\mu, \nu).
\]
\end{rmlist}
\begin{proof}
The proof of this fact is a generalisation of the method used in \cite{VAP98} to
prove the special case of the $\KL$ divergence.
Since $f(1) = 0$, we have the general property that
\[
f(x) = f(\max \{x, 1\}) + f(\min \{x, 1\}).
\]
Using this fact and the convexity of $f$ we get the general estimate 
\[
\begin{split}
D_f(\mu, \nu) & = \int f(\odq{\mu}{\nu}) \dd \nu \\
& = \int f(\max \{ \odq{\mu}{\nu}, 1 \}) \dd \nu
	+ \int f(\min \{ \odq{\mu}{\nu}, 1\}) \dd \nu \\
& \geq f (\int \max \{ \odq{\mu}{\nu}, 1 \} \dd \nu)
	+ f (\int \max \{ \odq{\mu}{\nu}, 1 \} \dd \nu).
\end{split}
\]
Now use that
\[
\begin{split}
\max \{ x, 1 \} & =  \frac{1 + x + |1 - x|}{2} \\
\min \{ x, 1 \} & =  \frac{1 + x - |1 - x|}{2}
\end{split}
\]
to complete the theorem.
\end{proof}
Recalling that always $\TV \leq 2$, the proposition rises the question as to when the function $f(1 + x) + f(1 - x)$ is monotonous on $x \in [0, 1]$.
The following lemma partially answers this.
\begin{rmlist}[Lemma]
Under the conditions of Lemma~\ref{rml:20.10}, the function $\phi(x) = f(1 + x) + f(1 - x)$ is strictly monotonous on $x \in [0, 1]$.
\end{rmlist}
\begin{proof}
The conditions imply that $\phi(0) = 0$, $\phi(x) > 0$ for  $x > 0$, and that $\phi$ is convex.
Let $0 \leq x_1 < x_2 \leq 1$.
For any $\tau \in ]0, 1[$, 
\[
(1 - \tau) \phi(0) + \tau \phi(x_2) > \phi((1 - \tau) 0 + \tau x_2)
\]
which obviously implies $\phi(x_2) > \tau \phi(x_2) >  \phi(\tau x_2)$ (since $\tau \in ]0, 1[$).
Now take $\tau = x_1 / x_2$ to get the result.
\end{proof}
As a corollary of Proposition~\ref{ftheorem}, we get the following well known estimates between $\TV$ and $\KL$
\begin{rmlist}[Corollary (Bretagnole--Huber and Furstemberg inequality)]
\label{bretagnolehuber}
\[
\TV(\mu , \nu) 
	\leq 2 \sqrt{1 - \exp \left(-\SH(\mu , \nu) \right)}
	\leq 2 \sqrt{\SH(\mu , \nu)}
\]
\end{rmlist}
Recall that $\SH(\mu , \nu) = \KL(\nu, \mu)$.
A further useful estimate concerns the Hellinger divergence
\begin{rmlist}[Corollary]
\label{hellingertheorem}
For the Hellinger divergence $\HE$, the estimate
\beq{equ:20.25}
\TV \leq \left \{
\begin{array}{ll}
	2 - 2 \left(1 - \sqrt{\HE}\right)^2 & \; \mbox{if } \HE < 1 \\
	2 & \; \mbox{otherwise}
\end{array} \right.
\eeq
holds.
\end{rmlist}
\begin{proof}
Theorem~\ref{ftheorem} gives the inequality
\beq{equ:20.30}
\HE \geq \left( \sqrt{ 1 + \frac{1}{2} \TV} - 1\right)^2 + \left(\sqrt{ 1 - \frac{1}{2} \TV} - 1\right)^2.
\eeq
The right hand side of Equation~\eqref{equ:20.30} is larger than $\left(\sqrt{1 - \frac{1}{2}\TV} - 1\right)^2$, whence
\[
\HE \geq \left(\sqrt{1 - \frac{1}{2} \TV} - 1 \right)^2,
\]
which, after solving for $\TV$, yields the result.
\end{proof}
\bibliographystyle{plain}
\bibliography{/home/broecker/TeX/Literatur}
\end{document}

%% file: main.macros.tex
\usepackage{epsfig,float,enumerate,ifthen,amsmath,amssymb,amsthm,delarray}
\newtheoremstyle{itlist}%
{}{}					
{\itshape}		
{}{}					
{.}{ }				
{\thmnumber{#2.}\thmnote{ \textit{#3}}} 
\newtheoremstyle{rmlist}%
{}{}					
{\normalfont}	
{}{}					
{.}{ }				
{\thmnumber{#2.}\thmnote{ \textit{#3}}} 
{\theoremstyle{itlist}
}
{\theoremstyle{rmlist}
\newtheorem{rmlist}[subsection]{}}

\newcommand{\R}{\mathbb R}

\newcommand{\vr}[1]{\left \langle #1 \right \rangle}  

\newcommand{\beq}[1]{\begin{equation}\label{#1}}
\newcommand{\eeq}{\end{equation}}
\newcommand{\beqn}{\begin{equation}\nonumber}
\newcommand{\dd}{\mathrm d}

\newcommand{\odq}[2]{\frac{\dd #1}{\dd #2}}
\newcommand{\TV}{\mathsf{TV}}
\newcommand{\HE}{\mathsf{HE}}

\newcommand{\KL}{\mathsf{KL}}
\newcommand{\PE}{\mathsf{PE}}
\newcommand{\SH}{\mathsf{SH}}
\ifthenelse{\value{the_style}=1}{%
}{\ifthenelse{\value{the_style}=2}{%
\hoffset-1in
\voffset-1in
\topmargin7mm
\headheight5mm
\headsep5mm
\textheight230mm
\oddsidemargin17mm
\textwidth125mm
\evensidemargin68mm
\usepackage{showlabels}
}{\ifthenelse{\value{the_style}=3}{%
\twocolumn
}{\typeout{No such style}}}}

%% file: main.bbl
\begin{thebibliography}{1}

\bibitem{csiszar67}
Imre Csiszar.
\newblock Information-type measues of difference of probability distributions
  and indirect observations.
\newblock {\em Studia Sci. Math. Hungar.}, 2:299--318, 1967.

\bibitem{doob94}
Joseph~L. Doob.
\newblock {\em Measure Theory}.
\newblock Springer, 1994.

\bibitem{VAP98}
Vladimir~N. Vapnik.
\newblock {\em Statistical Learning Theory}.
\newblock John Wiley \& Sons, Inc., New York, 1998.

\end{thebibliography}
